\documentclass[11pt, a4paper]{article}
\usepackage{amsfonts}
\usepackage{mathrsfs}
\usepackage{amssymb}
\usepackage{amsmath}
\begin{document}

\title{Approximation by Semigroups of Spherical Operators\thanks{The research was supported by
the National Natural Science Foundation of China (No. 60873206), the Natural Science Foundation of Zhejiang Province of
China (No. Y7080235) and the Innovation Foundation of Post-Graduates of Zhejiang Province of China (No.
YK2008066).}}

\author{Yuguang Wang \and Feilong Cao\thanks{Corresponding author: Feilong Cao,  E-mail: \tt feilongcao@gmail.com} }

 \date{}
\maketitle

\begin{center}
\footnotesize Department of Mathematics,
 China Jiliang University,

 Hangzhou 310018, Zhejiang Province,  P R China

\begin{abstract}
This paper discusses the approximation by 
a class of so called exponential-type multiplier operators. It is proved that such operators form a strongly continuous semigroup of contraction operators of  class ($\mathscr{C}_0$), from which the equivalence between approximation for these operators and $K$-functionals introduced by the operators is given. As examples, the constructed $r$-th Boolean of generalized spherical Abel-Poisson operator and $r$-th Boolean of generalized spherical Weierstrass operator denoted by $\oplus^r V_t^{\gamma}$ and $\oplus^r W_t^{\kappa}$ separately ($r$ is any positive integer, $0<\gamma,\kappa\leq1$ and $t>0$) satisfy that
    $\|\oplus^r V_t^{\gamma}f - f\|_{\mathcal{X}}\approx \omega^{r\gamma}(f,t^{1/\gamma})_{\mathcal{X}}$ and $\|\oplus^r W_t^{\kappa}f - f\|_{\mathcal{X}}\approx \omega^{2r\gamma}(f,t^{1/(2\kappa)})_{\mathcal{X}}$,
for all $f\in \mathcal{X}$, where $\mathcal{X}$ is a Banach space of continuous functions or $\mathcal{L}^p$-integrable functions ($1\leq p<\infty$) and $\|\cdot\|_{\mathcal{X}}$ is the norm on $\mathcal{X}$ and $\omega^s(f,t)_{\mathcal{X}}$ is the moduli of smoothness of degree $s>0$ for $f\in \mathcal{X}$.
The saturation order and saturation class of the regular exponential-type multiplier operators with positive kernels are also obtained. Moreover, it is proved that $\oplus^r V_t^{\gamma}$ and $\oplus^r W_t^{\kappa}$ have the same saturation class if $\gamma=2\kappa$.


\end{abstract}

\end{center}
\normalsize
\section{Introduction}
The study on spherical approximation started early in 1960's when Butzer, Berens and Pawelke (see \cite{Berens_Butzer_Pawelke1968}) studied the saturation properties of singular integrals on the sphere. Since 1980's, the approximation theory on the sphere has been developed by Nikol'ski$\breve{\i}$, Lizorkin et al. 
and since 1990's, Wang, Li and Dai et al. pursued the research of approximation theory on the sphere (see for example \cite{Wang_Li2006}). Some classical theorems in the case of one dimension, for example, Jackson-type theorem, were generalized to the sphere (see \cite{Riemenschneider_Wang1995}). The approximation tools then were mainly polynomial operators, say, spherical Jackson operators (see \cite{Li_Yang1991}) and de la Vall\'{e}e Poussin means (see \cite{Berens_Li1993}). In recent years, the study on spherical approximation has attracted more and more attention of researchers. There have been many interesting works in this field, such as \cite{Dai_Ditzian2005}-\cite{Ditzian2006}, \cite{Gia_Narcowich_Ward_Wendland2006}. In particular, we notice that as a non-polynomial operator on the sphere, the classical Abel-Poisson operator was studied by Dai and Ditzian \cite{Dai_Ditzian2005}, where the equivalence between approximation by Abel-Poisson operators on the sphere and the $K$-functional introduced by their infinitesimal generator was given.

This paper is mainly about the approximation by non-polynomial operators on the sphere, namely, a class of exponential-type operators denoted by $\{T_{p}^{\gamma}(t)|0\leq t<\infty\}$ with polynomial $p(x)$ and $0<\gamma \leq1$, in the form of
    $T_{p}^{\gamma}(t)f=\sum\limits_{k=0}^{\infty}e^{-(p(x))^{\gamma}t}Y_{k}f\;(f\in \mathcal{X})$,
where $Y_{k} f$ is the $k$-th term of the Laplace expansion of $f$ on the sphere. $T_{p}^{\gamma}(t)$ is called regular if the coefficient of the first item of $p(x)$ is positive, $p(0)=0$ and the degree of $p(x)$ is larger than $0$. These regular operators with positive kernels are proved to form a semigroup of class ($\mathscr{C}_{0}$).

The semigroups of operators were early studied by Hille and Phillips (see \cite{Hille_Phillips1957}). Later, in 1960's, Butzer and Berens studied the approximation properties of semigroups of operators (see \cite{Berens_Butzer_Pawelke1968}). 
By Ditzian and Ivanov's method (see \cite[Sec.5]{Ditzian_Ivanov1993}), we will prove that for a class of exponential-type multiplier operators that form a semigroup of class $(\mathscr{C}_{0})$ and for $r\in \mathbb{Z}_{+}$ there holds
\begin{equation}\label{eq12}
    \|(T(t)-I)^rf\|_{\mathcal{X}}\approx \inf_{g\in \mathcal{D}_{1}(\mathcal{A}^r)}(\|f-g\|_{\mathcal{X}}+\| \mathcal{A}^r g\|_{\mathcal{X}})
\end{equation}
if $\{T(s)|s>0\}$ possesses additive properties that $T(t)f\in \mathcal{D}_{1}(\mathcal{A})$ for all $t>0$ and $f\in \mathcal{X}$ and $\mathcal{A}$ satisfies the further following Bernstein-type inequality
\begin{equation}\label{eq77}
    t\|\mathcal{A}T(t)f\|_{\mathcal{X}}\leq N\|f\|_{\mathcal{X}}\quad(t>0,\;f\in\mathcal{X}),
\end{equation}
where $\mathcal{D}_{1}(\mathcal{A}^r)=\{g\in \mathcal {X} |~ \mathcal {A}^{r}g\in \mathcal {X}\}$ and $r$-th power of multiplier operator $\mathcal {A}$.

It will be proved that the regular $T_{p}^{\gamma}(t)$ with positive kernel also satisfies (\ref{eq77}). Hence (\ref{eq12}) holds for $T_{p}^{\gamma}(t)$. Noticing
$(I-T(t))^r f=f-\oplus^{r}T(t)f$,
here $\oplus^{r}T(t)$ is defined as the $r$-th Boolean of $T(t)$,
we shall obtain the equivalence between the approximation of $\oplus^{r}T_{p}^{\gamma}(t)$ and $K$-functional introduced by the multiplier operator, actually infinitesimal generator of $\big\{T_{p}^{\gamma}(t)|0\leq t<\infty\big\}$. In terms of convolution, we shall also obtain the saturation order and saturation class of Booleans of the regular $\oplus^{r}T_{p}^{\gamma}(t)$.
The generalized spherical Abel-Poisson operators $V_{t}^{\gamma}$ and generalized spherical Weierstrass operators $W_{t}^{\kappa}$ that were introduced by Bochner in \cite[P. 43-47]{Bochner1955} and \cite[P. 84]{Bochner1950} respectively, are actually the examples of regular exponential-type operators with positive kernels. Thus follows the equivalence between approximation of $\oplus^{r}V_{t}^{\gamma}$ or $\oplus^{r}W_{t}^{\kappa}$ and 
the moduli of smoothness on the sphere. It is also discussed that the saturation properties of $\oplus^{r}V_{t}^{\gamma}$ and $\oplus^{r}W_{t}^{\kappa}$.

The paper is structured as follows. In Section~2, some basic concepts will be introduced. Section~3 discusses the properties of some spherical function classes as well as $K$-functionals that are introduced by multiplier sequences. 
Section~4 consists of our main results and their proofs. In Section~5, the results of previous sections are applied to $\oplus^r V_t^{\gamma}$ and $\oplus^r W_t^{\kappa}$. Our main results are Theorem~4.2.3, Theorem~4.2.4, Theorem~5.3, Theorem~5.4 and Theorem~5.5.

\section{Preliminaries}
Denote by letters $C$, $C_i$ or $C(i)$ positive constants, where $i$ is either a positive integer, variable, function or space on which $C$ depends only. Their values may be different at different occurrences, even within one formula. The notation $a \approx b$ means that there exists a positive constant $C$ such that $C^{-1}b\leq a \leq Cb$. Denote by $f(t)=\mathcal{O}(t)$ which means there exists some constant $C$ independent of $t$ such that $|f(t)|\leq C |t|$, here $f(t)$ is a function with respect to $t$ and we write $f(t)=o(t)$ if $f(t)/t$ tends to zero as $t\rightarrow\infty$ or as $t\rightarrow t_{0}$ where $t_{0}$ is a real number. The collection of all positive integers are denoted by $\mathbb{Z}_{+}$. Let $\mathbb{S}^n$ be the unit sphere in $(n+1)$-dimensional Euclidean space $\mathbb{R}^{n+1}$. $x$ and $y$ are denoted as the points on $\mathbb{S}^{n}$ and $x\cdot y$ denotes the inner product in $ \mathbb{R}^{n+1}$. Denote by $d\omega_n(x)$ the elementary surface piece on $\mathbb{S}^n$ and by $d\omega(x)$ for convenience if there's no confusion. The volume of $\mathbb{S}^{n}$ is
$$
\left|\mathbb{S}^{n}\right|:=\int_{\mathbb{S}^{n}}d\omega(x)=\displaystyle\frac{2\pi^{n/2}}{\Gamma(n/2)}.
$$
Denote by $\mathcal{L}^p(\mathbb{S}^{n})$ the Banach space of $p$-th integrable
functions $f: \mathbb{S}^{n}\rightarrow\mathbb{C}$ ($\mathbb{C}$ is the collection of all complexes) with norm $
\| f\|_{\infty}:=\|f\|
_{\mathcal{L}^{\infty}(\mathbb{S}^{n})}:=\mbox{ess}\!\sup\limits_{x\in
\mathbb{S}^{n}}|f(x)|
$
and
$$
\|f\|_{\mathcal{L}^p}:=\|f\|_{\mathcal{L}^p(\mathbb{S}^{n})}:=\left\{\int_{\mathbb{S}^{n}}|f(x)|^pd\omega(x)\right\}^{1/p}<\infty
\quad
(1\leq p<\infty).
$$
Denote by $\mathcal{C}(\mathbb{S}^n)$ the Banach space consisting of all continuous functions $f: \mathbb{S}^n\rightarrow\mathbb{C}$ with norm
$
\|f\|_{\mathcal{C}}:=\max\limits_{x\in\mathbb{S}^{n}}|f(x)|.
$
Denote by $\mathcal{M}(\mathbb{S}^n)$ the collection of all finite regular Borel measures on $\mathbb{S}^n$ (the range is also in $\mathbb{C}$) and it is a Banach space with norm
$
\|\mu\|_{\mathcal{M}}:=\int_{\mathbb{S}^{n}}|d\mu(x)|.
$
$\mathcal{L}^p(\mathbb{S}^n)\;(1\leq p\leq\infty)$, $\mathcal{C}(\mathbb{S}^n)$ and $\mathcal{M}(\mathbb{S}^n)$ may be replaced by $\mathcal{L}^p$, $\mathcal{C}$ and $\mathcal{M}$ for convenience. Denote by $\mathcal{X}$ either $\mathcal{L}^p(\mathbb{S}^n)\;(1\leq p<\infty)$ or $\mathcal{C}(\mathbb{S}^n)$. The dual space of $\mathcal{X}$, the collection of all bounded linear functionals on $\mathcal{X}$, is denoted by $\mathcal{X}^*$.

A function $f$ on $\mathbb{R}^{n+1}$ is called harmonic if $\Delta f=0$ where
$\Delta=\frac{\partial}{\partial x_{1}^2}+\frac{\partial}{\partial x_{2}^2}+\dots+\frac{\partial}{\partial x_{n+1}^2}$
is the classical Laplace operator on $\mathbb{R}^{n+1}$. The collection of all homogeneous and harmonic polynomials on $\mathbb{R}^{n+1}$ is denoted by $\mathcal{A}_k^n$. And the collection of restrictions on $\mathbb{S}^{n}$ of all functions in $\mathcal{A}_k^n$ is denoted by $\mathcal{H}_{k}^{n}$.
Denote by $\Pi_{k}^n$ the collection of restrictions on $\mathbb{S}^{n}$ of all polynomials on $\mathbb{R}^{n+1}$ which is dense in $\mathcal{X}$
and any polynomial restricted on $\mathbb{S}^{n}$ with degree $k\in \mathbb{Z}_{+}$ is in $\mbox{\rm span}\{\mathcal{H}_{j}^{n}| 0\leq j\leq k\}$, the linear combination of $\mathcal{H}_{1}^{n}$, $\mathcal{H}_{2}^{n}$, $\dots$, $\mathcal{H}_{k}^{n}$.\\
\textbf{Definition~2.1}~~\emph{The $r$-th {\rm($r\in \mathbb{Z}_{+}$)} Boolean of an operator $T$ on $\mathcal{X}$ (an operator from $\mathcal {X}$ to $\mathcal {X}$) is defined as
\begin{equation*}\label{eq25}
    \oplus^{r}T:=I-(I-T)^r=-\sum_{i=1}^{r}(-1)^{i}{r\choose i}T^{i},
\end{equation*}
where $\displaystyle{r\choose k}:=\frac{r!}{k!(r-k)!}$ and $T^{0}:=I$.}\\[0.05 cm]
The projection on $\mathcal{H}_{k}^{n}$ of $f\in \mathcal{L}^{1}(\mathbb{S}^n)$ is defined by (see \cite[Chap.1]{Berens_Butzer_Pawelke1968} and \cite[Chap.1]{Wang_Li2006})
\begin{eqnarray*}
Y_k(f)(x):=\frac{\Gamma(\lambda)(k+\lambda)}{2\pi^{\lambda+1}}\int_{\mathbb{S}^{n}}P_k^{\lambda}(x\cdot
y)f(y)\:d\omega_{n}(y)
\end{eqnarray*}
and for $\mu\in\mathcal{M}(\mathbb{S}^n)$,
 $Y_k(d\mu)(x):=\frac{\Gamma(\lambda)(k+\lambda)}{2\pi^{\lambda+1}}\displaystyle\int_{\mathbb{S}^{n}}P_k^{\lambda}\;(x\cdot
y)\:d\mu(y)$,
where $2\lambda=n-2$, and $P_k^{\nu}(t)$, $|t|\leq1$, $k=0,1,2,\dots$, $\nu>-1/2$ is the ultraspherical polynomial (Gegenbauer polynomial) of degree $k$ with $\nu$ and is generated by
    $(1-2 t r+r^2)^{-\nu}=\sum\limits_{k=0}^{\infty}P_k^{\nu}(t)r^k\; (0\leq r<1)$.
$\{P_k^{\nu}(t)\}_{k=0}^{\infty}$ forms an orthogonal system with the weight $(1-t^2)^{\nu-1/2}$, that is, for $\nu>-1/2,\; \nu\neq 0$ (see \cite[P. 81]{Szego2003}),
\[\begin{array}{ll}
  & \displaystyle\int_{-1}^{1}P_k^{\nu}(t)P_j^{\nu}(t)(1-t^2)^{\nu-1/2}\:dt \\[0.2 cm]
  =& \displaystyle\int_{0}^{\pi} P_k^{\nu}(\cos\theta)P_j^{\nu}(\cos\theta)(\sin\theta)^{2\nu}\:d\theta
  =\left\{\begin{array}{ll}
  (c(k,\nu))^{-1} & (k=j), \\
  0 &(k\neq j).
  \end{array}\right.
\end{array}\]
\vspace{-2 cm}\begin{equation}\label{eq33}\end{equation}\\
where
    $c(k,\nu):=(2^{2\nu-1}(\Gamma(\nu))^2(k+\nu)\Gamma(k+1))/(\pi \Gamma(k+2\nu))$.
Now let $\nu=\lambda :=(n-2)/2>0$, then (see \cite[P. 171]{Szego2003})
\begin{equation}\label{eq78}
    \left|P_{k}^{\lambda}(t)\right|=\mathcal {O}(k^{2\lambda-1}).
\end{equation}


A function $f\in \mathcal{X}$ is called a zonal function with $x_0$ on $\mathbb{S}^n$ if for some fixed $x_0\in \mathbb{S}^n$, $f(x_0\cdot y)$ is a constant when $x_0\cdot y$ is unchanged. The collection of all zonal functions with $x_0$ in $\mathcal{L}^{p}$ or $\mathcal{C}$ is denoted by $\mathcal{L}^p_{\lambda}(\mathbb{S}^{n},x_0)$ $(1\leq p<\infty)$, $\mathcal{C}_{\lambda}(\mathbb{S}^{n},x_0)$ ($\mathcal{L}^p_{\lambda}$, $\mathcal{C}_{\lambda}$ for convenience if there's no confusion).
$\mathcal{L}^p_{\lambda}(\mathbb{S}^{n},x_0)$ $(1\leq p<\infty)$ with norm
\begin{equation}\label{eq74}
    \|\varphi\|_{\mathcal{L}^p_{\lambda}}:=\left\{\int_{\mathbb{S}^n}|\varphi(x\cdot y)|^pd\omega(y)\right\}^{1/p}=\left\{\left|\mathbb{S}^{n-1}\right|\int_{0}^{\pi}|\varphi(\cos\theta)|^p(\sin\theta)^{2\lambda}d\theta\right\}^{1/p},
\end{equation}
$\mathcal{C}_{\lambda}(\mathbb{S}^{n},x_0)$ with norm
    $\|\varphi\|_{\mathcal{C}_{\lambda}}:=\sup\limits_{0\leq\theta\leq\pi} |\varphi(\cos\theta)|$,
and $\mathcal{M}_{\lambda}(\mathbb{S}^{n},x_0)$ with norm
    $\|\mu\|_{\mathcal{M}_{\lambda}}:=\left|\mathbb{S}^{n-1}\right|\displaystyle\int_{0}^{\pi}|d\mu^{*}(\theta)|$,
where $\mu^{*}$ is the corresponding function in $\mathcal{M}[0,\pi]$ of the measure $\mu\in\mathbb{S}^n$ (actually there is a bijection between $\mathcal{M}_{\lambda}(\mathbb{S}^{n},x_0)$ and some subset of $\mathcal{M}[0,\pi]$, see \cite[P. 250-252]{Dunkl1966} and \cite[P. 204 and P. 209]{Berens_Butzer_Pawelke1968}), are all Banach spaces.

For $f\in \mathcal{L}^1(\mathbb{S}^n)$ and $\varphi\in \mathcal{L}^1_{\lambda}(\mathbb{S}^n)$, the convolution of $f$ and the zonal function $\varphi$ is defined by
\begin{equation}\label{eq35}
    (f*\varphi)(x):=\int_{\mathbb{S}^n}f(y)\varphi(x\cdot y)\:d\omega(y)\quad(x\in \mathbb{S}^n).
\end{equation}

The convolution of $\psi\in \mathcal{L}^1_{\lambda}(\mathbb{S}^n)$ and $\mu\in \mathcal{M}(\mathbb{S}^{n})$ is defined by
    $(\psi*d\mu)(x):=\displaystyle\int_{\mathbb{S}^{n}}\psi(x\cdot y)d\mu(y).$
The convolution of $f\in \mathcal{L}^1(\mathbb{S}^n)$ and the zonal measure $\mu\in \mathcal{M}_{\lambda}(\mathbb{S}^n)$ with $x_0$ is defined by
\begin{equation}\label{eq36}
    (f*d\mu)(x):=\int_{\mathbb{S}^n}f(y)\:d\varphi_{x}\mu(y)\quad(x\in \mathbb{S}^n),
\end{equation}
where $\varphi_{x}\mu(E):=\mu(\rho E)$, $\rho x=x_0$, for all measurable subsets $E\subset \mathbb{S}^n$ (
Please refer to \cite[Chapter 1]{Berens_Butzer_Pawelke1968}, \cite[Chapter 1]{Wang_Li2006} and \cite{Dunkl1966} for further details on convolution).\\
\textbf{Remark~2.2} 
~~\emph{In this paper, we follow the definition of convolution in \cite{Wang_Li2006} 
and Young's inequality still holds on the sphere.}\\
\textbf{Definition~2.3} (see \cite[P. 254]{Berens_Butzer_Pawelke1968})~~\emph{Let two function spaces $\mathcal{Y}$ and $\mathcal{Z}$ either be $\mathcal{C}(\mathbb{S}^{n})$, $\mathcal{L}^p(\mathbb{S}^{n})$ {\rm($1\leq p<\infty$)} or $\mathcal{M}(\mathbb{S}^{n})$. A sequence $\{a_k\in \mathbb{C}| \:k=0,1,2,\dots\}$ is called a multiplier sequence from $\mathcal {Y}$ to $\mathcal {Z}$ if for each $f\in \mathcal{Y}$, there exists $g\in \mathcal{Z}$ whose Laplace expansion is as follows
    $Y_k g=\lambda/(k+\lambda)\:a_{k} Y_k f\;(k=0,1,2,\dots)$.
The collection of all multiplier sequences from $\mathcal {Y}$ to $\mathcal {Z}$ is denoted by $(\mathcal{Y},\mathcal{Z})$. For $\mathcal{Y}=\mathcal{Z}$, $\{a_k\}_{k=0}^{\infty}$ is called a multiplier sequence on $\mathcal{Y}$.}\\
\textbf{Remark~2.4}~~\emph{By {\rm\cite[P. 222-P. 231]{Kaczmarz_Steinhaus1935}}
\begin{equation*}
    (\mathcal{M},\mathcal{M})=(\mathcal{C},\mathcal{C})=(\mathcal{L}^{1},\mathcal{L}^{1})\subset(\mathcal{L}^{p},\mathcal{L}^{p})\subset(\mathcal{L}^{2},\mathcal{L}^{2})\quad(1<p<\infty).
\end{equation*}}
\textbf{Definition~2.5}~~\emph{An operator $T$ on $\mathcal {X}$ is called a multiplier operator with a sequence $\{a_k\}_{k=0}^{\infty}$ on $\mathcal {X}$ if for each $f\in \mathcal{X}$, $Tf\in \mathcal {X}$ and
    $Y_k(Tf)=a_k Y_k(f)\;(k=0,1,2,\dots)$.
}
\emph{The operator $T^{\alpha}$ ($\alpha>0$) defined by
    $T^{\alpha} f \sim \sum_{k=0}^{\infty}-(a(k))^{\alpha} Y_k(f)$,
where ``\;$\sim$'' is in the sense of distribution (see \cite[P. 323-325]{Ditzian1998} and \cite[Section.1]{Grafakos2005}), is called the fractional differential operator if~ $T^{\alpha} f\in \mathcal {X}$. Denote the domain of $T^{\alpha}$ by
    $\mathcal {D}_{1}(T^{\alpha})=\left\{f\in \mathcal {X}\big|~ T^{\alpha}f\in \mathcal {X}\right\}$.
}\\
\hspace{1.3 em}Denote by $\|T\|_{\mathcal{X}}$ the norm of an operator $T$ on $\mathcal{X}$. Then the collection of endomorphisms of $\mathcal{X}$ denoted by $\mathscr{E}(\mathcal{X})$ is a Banach algebra with norm $\|T\|_{\mathcal{X}}$ (see \cite[P. 7]{Buter_Berens1967}).\\
\textbf{Definition~2.6} (see \cite[P. 7-8]{Buter_Berens1967})~~\emph{If~ $T(t)$ is an operator function on the non-negative real axis ~$0\leq t<\infty$ to the Banach algebra $\mathscr{E}(\mathcal{X})$, in the following conditions, if {\rm(\ref{eq15})} is satisfied, $\{T(t)|\:0\leq t<\infty\}$ is called one-parameter semigroup of operators in $\mathscr{E}(\mathcal{X})$ and it is said to be of class {\rm($\mathscr{C}_0$)} if it satisfies the further property {\rm(\ref{eq16})},
\begin{eqnarray}
 T(t_1+t_2) &=& T(t_1)\circ T(t_2)\;(t_1,t_2\geq0),\quad T(0) = I,\label{eq15}\\
 s\!-\!\!\!\lim_{t\rightarrow0+}T(t) &=& I,\label{eq16}
\end{eqnarray}
where $A\circ B$ means the composition of operators $A$ and $B$, $I$ is the identity operator on $\mathcal{X}$ and $s\!-\!\!\!\lim\limits_{t\rightarrow0+}f_{t}=f$ denotes the strongly convergence which means $\|f_{t}-f\|_{\mathcal {X}}$ tends to zero as $t\rightarrow0+$.
$T(t)$ is called to have contraction if it satisfies
\begin{equation}\label{eq81}
    \|T(t)f\|_{\mathcal{X}}\leq\|f\|_{\mathcal{X}}\quad(f\in \mathcal{X}).
\end{equation}}
\textbf{Definition~2.7} (see \cite[P. 11]{Buter_Berens1967})~~\emph{The infinitesimal generator $\mathcal{A}$ of the semigroup $\{T(t)|\,0\leq t<\infty\}$ is defined by
    $\mathcal{A}:=s\!-\!\!\!\lim\limits_{t\rightarrow0+}\displaystyle\frac{T(t)f - f}{t}$,
whenever the limit exists; the domain of $\mathcal{A}$ is, in symbols $\mathcal{D}(\mathcal{A})$, being the set of elements $f\in \mathcal{X}$ for which the limit exists; for $r=0,1,2,\dots$, the $r$-th power of $\mathcal{A}$ denoted by $\mathcal{A}^{r}$ is defined inductively by the relations $\mathcal{A}^0=I$, $\mathcal{A}^1=\mathcal{A}$, and\\[-0.3 cm]
\parbox{0.9\textwidth}{\begin{eqnarray*}
  && \mathcal{D}(\mathcal{A}^{r}) := \Big\{f|\:f\in \mathcal{D}(\mathcal{A}^{r-1})\; \mbox{and}\; \mathcal{A}^{r-1}f\in \mathcal{D}(\mathcal{A})\Big\}, \nonumber\\
  && \mathcal{A}^{r}f := \mathcal{A}(\mathcal{A}^{r-1}f)=s\!-\!\!\!\lim_{t\rightarrow0+}\frac{T(t)-I}{t}\mathcal{A}^{r-1}f\quad(f\in \mathcal{D}(\mathcal{A}^{r})).
\end{eqnarray*}}\hfill
\parbox{0.1 \textwidth}{
\begin{equation}\label{eq18}
\end{equation}
}}
For $r\in \mathbb{Z}_{+}$, $\mathcal{D}(\mathcal{A}^{r})$ is a linear subspace and $\mathcal{A}^{r}$ is a linear operator.


If an operator $T$ in $\mathscr{E}(\mathcal{X})$ can be expressed in the form of convolution (\ref{eq35}) or (\ref{eq36}), then $\varphi$, $\psi\in \mathcal{L}_{\lambda}^1$ or $\mu\in \mathcal{M}_{\lambda}$ there is called the kernel of $T$.
The Ces\`{a}ro mean of $f\in\mathcal{X}$ denoted by $\sigma_k^{\alpha}(f)$ is defined by (see for instance \cite[P. 49]{Wang_Li2006})
    $\sigma_k^{\alpha}(f)=(1/A_{k}^{\alpha})\sum_{j=0}^{k}A_{k-j}^{\alpha}Y_{j}f$,
where $\alpha$ is a complex whose real part is not less than $-1$, $k\in \mathbb{Z}_{+}$ and
 $A_{k}^{\alpha}={k+\alpha\choose\alpha}=\Gamma(k+\alpha+1)/(\Gamma(\alpha+1)\Gamma(k+1))\;(k\in \mathbb{Z}_{+})$,
is the generalized combination number.
For $\alpha>\lambda=(n-2)/2$, there holds (see for instance Theorem~2.3.10 in \cite[P. 54-55]{Wang_Li2006})
\begin{equation}\label{eq46}
    \|\sigma_k^{\alpha}(f)\|_{\mathcal{X}}\leq
    C(n,\alpha,\mathcal{X})\|f\|_{\mathcal{X}}\quad(k\in\mathbb{Z}_{+},\;f\in \mathcal{X}).
\end{equation}
For any sequence $\{\mu_k\}_{k=0}^{\infty}$,
denote $\delta\mu_{k}=\mu_{k}-\mu_{k+1}$ $(k=0,1,\dots)$ and $\delta^{i+1}\mu_{k}=\delta(\delta^{i}\mu_k)$ $(i=1,2,\dots)$.

The definitions of moduli of smoothness on the sphere are given as follows (see for instance \cite[P. 56- P.57, P. 183-184]{Wang_Li2006}). The translation operator on $\mathcal{L}^1(\mathbb{S}^n)$ is defined by
 $S_\theta(f)(x):=\big(|\mathbb{S}^{n-2}|\sin^{n-1}\theta\big)^{-1}\displaystyle\int_{x\cdot y=\cos \theta}f(y)\:d\omega_{n-1}(y)\;(0<\theta\leq\pi)$.
Let $\alpha>0$, $\theta>0$. The multiplier operator on $\mathcal{X}$ is called the finite difference of degree $\alpha$ with step $\theta$, defined by
    $\Delta_{\theta}^{\alpha}:=(I-S_{\theta})^{\frac{\alpha}{2}}=\sum\limits_{k=0}^{\infty}(-1)^k {\frac{\alpha}{2}\choose k}(S_{\theta})^k$,
where
    ${\frac{\alpha}{2}\choose k}:=\frac{1}{k!}\frac{\alpha}{2}\left(\frac{\alpha}{2}-1\right)\cdots\left(\frac{\alpha}{2}-k+1\right)$.\\
\textbf{Definition~2.8}~~\emph{Let $f\in \mathcal{X}$, $\alpha>0$. The moduli of smoothness of degree $\alpha$ of $f$ is defined as
    $\omega^{\alpha}(f,t)_{\mathcal{X}}:=\sup\Big\{\|\Delta_{\theta}^{\alpha} f\|_{\mathcal{X}}: 0<\theta\leq t\Big\}\;(0<t\leq\pi)$.
}\\
\textbf{Definition~2.9} (see for instance \cite[P. 323-325]{Ditzian1998})~~\emph{Let $f\in \mathcal{X}$, $t>0$. The $K$-functional introduced by multiplier operator $\mathcal{A}$ with multiplier sequence $\{a_k\}_{k=0}^{\infty}$ is defined as
    $K_{\mathcal{A}}(f,t)_{\mathcal{X}}:=\inf\limits_{g\in \mathcal{D}_{1}(\mathcal{A})}\Big\{\|f-g\|_{\mathcal{X}}+t\|\mathcal{A}g\|_{\mathcal{X}}\Big\}$,
} here
\begin{eqnarray}\label{eq80}
    \mathcal{D}_{1}(\mathcal{A})=\bigg\{f\in \mathcal{X}|
    \mbox{\rm there exists}~~g\in \mathcal{X}~~\mbox{\rm such that}~~
    a(k)Y_{k}f=Y_{k}g,\;k=0,1,2,\dots\bigg\}.
\end{eqnarray}

Particularly, for the $K$-functional introduced by $(\alpha/2)$-th Laplace-Beltrami operator $D^{\alpha/2}$ with the multiplier sequence $\Big\{\big(-k(k+2\lambda)\big)^{\alpha/2}\Big\}_{k=0}^{\infty}$ ($\alpha>0$), there holds
\begin{equation}\label{eq45}
    K_{D^{\alpha/2}}(f,t^{\alpha})_{\mathcal{X}}\approx \omega^{\alpha}(f,t)_{\mathcal{X}},
\end{equation}
which was finally proved by Riemenschneider and Wang (see \cite{Riemenschneider_Wang1995}). One might also define $K$-functional introduced by infinitesimal generator as follows.\\
\textbf{Definition~2.10} (see \cite[Section~3.4]{Buter_Berens1967})~~\emph{Suppose $\{T(t)|0\leq t<\infty\}$ a semigroup of operator of class {\rm $(\mathscr{C}_{0})$} in {\rm$\mathscr{E}(\mathcal {X})$} and let $\mathcal {A}$ be its infinitesimal generator. For $f\in \mathcal {X}$, the $r$-th $K$-functional introduced by $\mathcal {A}$ is defined by
    $K_{\mathcal{A}^{r}}^{*}(f,t)_{\mathcal{X}}:=\inf\limits_{g\in \mathcal{D}(\mathcal{A}^{r})}\Big\{\|f-g\|_{\mathcal{X}}+t\|\mathcal{A}^{r}g\|_{\mathcal{X}}\Big\}$,
here $\mathcal{D}(\mathcal{A}^{r})$ is defined by Definition~2.7 and $\mathcal{A}^{r}$ denotes the $r$-th power of infinitesimal generator $\mathcal{A}$.}

Finally, it is worth mentioning here the concept of saturation for operators on $\mathcal {X}$ (see \cite[P. 217]{Berens_Butzer_Pawelke1968}), which was first proposed by Favard in \cite{Favard1949}.\\
\textbf{Definition~2.11}~~\emph{ Let $\varphi(\rho)$ be a positive function with respect to $\rho$, $0<\rho<\infty$, tending monotonely to zero as $\rho\rightarrow\infty$. For a sequence of operators $\{I_\rho\}_{\rho>0}$ if there exists $\mathcal {K}\subseteqq \mathcal{X}$ such that
\vspace{0.1 cm}\\\parbox{6 cm}{\begin{tabular}{ll}
(i) & If\; $\| I_\rho(f)-f\|_{p} =o(\varphi(\rho))$, then $I_{\rho}f=f$ for all $\rho>0$; \\
(ii) & $\| I_\rho(f)-f\|_{p} = \mathcal{O}(\varphi(\rho))$ if and only if
$f\in\mathcal {K}$;
\end{tabular}}\parbox{1 cm}{}\vspace{0.2 cm}\\
then $I_\rho$ is said to be saturated on $\mathcal {X}$ with order $\mathcal{O}(\varphi(\rho))$ and $\mathcal {K}$ is called its saturation class.}
\section{Classes and $K$-Functionals Introduced by Multiplier Sequences on the Sphere}
Let $\psi(x)$ be a function from $\mathbb{R}$ to $\mathbb{C}$, define 
    $\mathcal{H}\Big(\big\{\psi(k)\big\}_{k=0}^{\infty};\mathcal{X}\Big):=\Big\{f\in 
    \mathcal{X}\big|~\mbox{there exists}~g\in\mathcal{X}~\mbox{such that}
    ~\psi(k)Y_{k}f=Y_{k}g~\mbox{for}~k=0,1,\dots\Big\}$,
and denoted by $\mathcal{H}\big(\psi(k);\mathcal{X}\big)$ for convenience.\\
\textbf{Theorem~3.1}~~\emph{Suppose that $\psi_{0}(x)$ and $\varphi_{0}(x)$ are functions from $[0,+\infty)$ to $\mathbb{C}$ and there exist $v_1$, $v_2\in \mathbb{R}$ such that $\psi(x)=e^{iv_{1}\pi}\psi_{0}(x)$ and $\varphi(x)=e^{iv_{2}\pi}\varphi_{0}(x)$ are both real valued functions. And
 $0<\lim\limits_{x\rightarrow+\infty}\big(\psi(x)/\varphi(x)\big)=c_{0}<+\infty$
and $\psi(0)=\varphi(0)=0$, setting
\[g(t):=\left\{\begin{array}{ll}
\displaystyle\frac{\psi(t^{-1})}{\varphi(t^{-1})}, & 0<t<+\infty,\\
c_{0},& t=0,
\end{array}\right.\]
if $g(t),\;(g(t))^{-1}\in \mathcal{C}^{2\lambda+2}[0,+\infty)$ \mbox{\rm(}$\mathcal{C}^{2\lambda+2}[0,+\infty)$ is the collection of real functions on $[0,+\infty)$ that are $(2\lambda+2)$ times continuously differentiable\mbox{\rm)}, then for $-\infty<s<+\infty$, there holds
\begin{equation*}
    \mathcal{H}\big((\psi_{0}(k))^s;\mathcal{X}\big)=\mathcal{H}\big((\varphi_{0}(k))^s;\mathcal{X}\big).
\end{equation*}}
\hspace{1.3 em}\textbf{Proof.} We first prove that for any $s\in (-\infty,+\infty)$,
\begin{equation*}
    \left\{C_{k}^{s}=\frac{k+\lambda}{\lambda}\bigg(\frac{\psi(k)}{\varphi(k)}\bigg)^s, k=1,2,\dots;C_{0}^{s}=\varphi(0)\right\}
\end{equation*}
belongs to $(\mathcal{M},\mathcal{M})$.
In fact, for any $k\in \mathbb{Z}_{+}$, $g_1(t)=(g(t))^s\in \mathcal{C}^{2\lambda+2}[0,+\infty)$ allows us to use Taylor's formula for $g_1(t)$ on $\left[0,1/k\right]$ at $t=0$, that is, there exists $0<\xi_{k}<1/k$ such that
\begin{eqnarray}\label{eq29}
  \left(\frac{\psi(k)}{\varphi(k)}\right)^{s} &=& g_{1}\left(\frac{1}{k}\right)=g_{1}(0)+g_{1}^{(1)}(0)\frac{1}{k}+\cdots+\frac{g_{1}^{(2\lambda+1)}(0)}{(2\lambda+1)!}\left(\frac{1}{k}\right)^{2\lambda+1}\nonumber\\
  &&+\frac{g_{1}^{(2\lambda+2)}(\xi_{k})}{(2\lambda+2)!}\left(\frac{1}{k}\right)^{2\lambda+2}.
\end{eqnarray}
We deduce from the assumption that $g_1^{(i)}(0)$, $i=0,1,\dots,2\lambda+1$, are constants depending only on $\psi$, $\varphi$, $s$ and $n$, and
\begin{equation}\label{eq30}
\left|g_{1}^{(2\lambda+2)}(\xi_{k})\right|\leq C(\varphi, \psi, s, n).
\end{equation}
Multiply (\ref{eq29}) by $(n+\lambda)/\lambda$, then according to Definition~2.3, one can verify that the sequence consisting of the first term $\big(g_{1}(0)(k+\lambda)\big)/\lambda=\big(c_{0}(k+\lambda)\big)/\lambda$ $(k=1,2,\dots)$ belongs to $(\mathcal{M},\mathcal{M})$. \cite[I, P. 202-203]{Askey_Wainger1965} proved that $(k+\lambda)/k^{\alpha}$ $(\alpha>0)$ are Gegenbauer-Stieltjes-coefficients of some measure in $\mathcal{M}$. For the last term of (\ref{eq29}), we estimate the following series that
\begin{eqnarray*}
  \left|\frac{1}{|\mathbb{S}^n|}\sum_{k=1}^{\infty} \frac{g_{1}^{(2\lambda+2)}(\xi_{k})}{(2\lambda+2)!}\left(\frac{1}{k}\right)^{2\lambda+2}\frac{k+\lambda}{\lambda} P_{k}^{\lambda}(\cos\theta)(\sin\theta)^{2\lambda}\right|&\leq& C(\varphi, \psi, s, n)\sum_{k=0}^{\infty}\frac{1}{k^{2}}<\infty,
\end{eqnarray*}
here the inequality is due to (\ref{eq30}) and (\ref{eq78}). Thus, there exits $\mu_{1}\in \mathcal{M}_{\lambda}(\mathbb{S}^n)$ such that
\begin{equation*}
    d\mu_{1}^{*}(\theta)=\frac{1}{|\mathbb{S}^n|}\sum_{k=1}^{\infty} \frac{g_{1}^{(2\lambda+2)}(\xi_{k})}{(2\lambda+2)!}\left(\frac{1}{k}\right)^{2\lambda+2}\frac{k+\lambda}{\lambda} P_{k}^{\lambda}(\cos\theta)(\sin\theta)^{2\lambda}d\theta.
\end{equation*}
It follows that the Gegenbauer-Stieltjes-coefficients of $\mu_{1}$ are
\begin{eqnarray*}
    \check{\mu}_{1}(j)&=&|\mathbb{S}^n|c(j,\lambda)\int_{0}^{\pi}P_{j}^{\lambda}(\cos\theta)d\mu_{1}^{*}(\theta)\\
    &=&\frac{j+\lambda}{\lambda}\frac{g_{1}^{(2\lambda+2)}(\xi_{j})}{(2\lambda+2)!}\left(\frac{1}{j}\right)^{2\lambda+2}\quad(j=1,2,\dots).
\end{eqnarray*}

Lemma~5.3.1 in \cite[P. 255]{Berens_Butzer_Pawelke1968} tells us that a sequence is Gegenbauer-Stieltjes-coefficients of some zonal measure on $\mathbb{S}^{n}$ if and only if it belongs to $(\mathcal{M},\mathcal{M})$, hence,
\begin{equation*}
 \left\{\frac{k+\lambda}{\lambda}\frac{g_{1}^{(2\lambda+2)}(\xi_{k})}{(2\lambda+2)!}\left(\frac{1}{k}\right)^{2\lambda+2}\right\}_{k=1}^{\infty}
\end{equation*}
are Gegenbauer-Stieltjes-coefficients of $\mu_{1}$, which implies that
\begin{equation*}
    \left\{C_{k}^{s}=\frac{k+\lambda}{\lambda}\bigg(\frac{\psi(k)}{\varphi(k)}\bigg)^s,\: k=1,2,\dots;C_{0}^{s}=0\right\}
\end{equation*}
belongs to $(\mathcal{M},\mathcal{M})$.
By Remark~2.4, we obtain
\begin{eqnarray}\label{eq31}
    \big\{C_{k}^{s}\big\}_{k=0}^{\infty}\in(\mathcal{M},\mathcal{M})&=&(\mathcal{C},\mathcal{C})
    \subset(\mathcal{L}^{p},\mathcal{L}^{p})
    \;(1\leq p<\infty).
\end{eqnarray}
Now we prove $\mathcal{H}\big((\varphi(k))^s;\mathcal{X}\big)=\mathcal{H}\big((\psi(k))^s;\mathcal{X}\big)$. We will just take account of the case of $\mathcal{X}=\mathcal{L}^{p}(\mathbb{S}^n)\;(1\leq p<\infty)$ and the proof of the case of $\mathcal{X}=\mathcal{C}(\mathbb{S}^n)$ is analogous.
For $f\in \mathcal{H}\big((\psi(k))^s;\mathcal{L}^{p}(\mathbb{S}^n)\big)$ ($1\leq p<\infty$, $s\in \mathbb{R}$), there exists $g_{1}\in \mathcal{L}^{p}(\mathbb{S}^n)$ such that
 $(\psi(k))^s\: Y_{k}f=Y_{k} \:g_{1}\;(k=0,1,2,\dots)$.
Thus,
\begin{eqnarray*}
    (\varphi(k))^s\: Y_{k}f 
    &=& \frac{\lambda}{k+\lambda}C_{k}^{-\!s}\:Y_{k}\:g_{1}\quad(k=1,2,\dots).
\end{eqnarray*}
It follows from (\ref{eq31}) that $C_{k}^{-\!s}\in (\mathcal{L}^{p},\mathcal{L}^{p})$. So, there exists $g_{2}\in \mathcal{L}^{p}(\mathbb{S}^n)$ such that
\begin{equation*}
    \frac{\lambda}{k+\lambda}C_{k}^{-\!s}\:Y_{k}\:g_{1}=Y_{k}\:g_{2}\quad(k=0,1,2,\dots),
\end{equation*}
that is,
    $(\varphi(k))^{s}\: Y_{k}f=Y_{k}\:g_{2}\;(k=1,2,\dots)$,
in addition, $(\varphi(0))^{s}Y_{k}f=0=Y_{0}\:g_{2}$, so $f\in \mathcal{H}\big((\varphi(k))^s;\mathcal{X}\big)$.
Thus,
 $\mathcal{H}\big((\psi(k))^s;\mathcal{X}\big)\subset\mathcal{H}\big((\varphi(k))^s;\mathcal{X}\big)$.
One will similarly obtain  that
    $\mathcal{H}\big((\varphi(k))^{s};\mathcal{X}\big)\subset\mathcal{H}\big((\psi(k))^{s};\mathcal{X}\big).$
Hence,
    $\mathcal{H}\big((\varphi_{0}(k))^{s};\mathcal{X}\big)=\mathcal{H}\big((\varphi(k))^{s};\mathcal{X}\big)
    =\mathcal{H}\big((\psi(k))^{s};\mathcal{X}\big)=\mathcal{H}\big((\psi_{0}(k))^{s};\mathcal{X}\big)$.
This completes the proof of Theorem~3.1.\quad$\Box$\\
\textbf{Remark~3.2}~~\emph{Define\\
\[\mathcal{H}_{1}\Big(\big\{\psi(k)\big\}_{k=0}^{\infty};\mathcal{X}\Big):=\left\{
\begin{array}{ll}
    \Big\{f\in &\!\!\!\! \mathcal{X}\big|~\mbox{there exists}~g\in \mathcal{L}^{p}(\mathbb{S}^{n})~\mbox{such that}~\psi(k)Y_{k}f=Y_{k}g~\\ &\mbox{for}~k=0,1,\dots\Big\}\;\left(\mathcal{X}=\mathcal{L}^{p}(\mathbb{S}^{n}),\:1<p<\infty\right),\\[0.2 cm]
    \Big\{f\in &\!\!\!\! \mathcal{X}\big|~\mbox{there exists}~\mu\in \mathcal{M}(\mathbb{S}^{n})~\mbox{such that}\\
    &~\psi(k)Y_{k}f=Y_{k}(d\mu)~\mbox{for}~k=0,1,\dots\Big\}\;\left(\mathcal{X}=\mathcal{L}^{1}(\mathbb{S}^{n})\right),\\[0.2 cm]
    \Big\{f\in &\!\!\!\! \mathcal{X}\big|~\mbox{there exists}~g\in \mathcal{L}^{\infty}(\mathbb{S}^{n})~\mbox{such that}\\
    &~\psi(k)Y_{k}f=Y_{k}g~\mbox{for}~k=0,1,\dots\Big\}\;\left(\mathcal{X}=\mathcal{C}(\mathbb{S}^{n})\right).
\end{array}\right.\]\\
Suppose $\varphi_{0}(x)$ and $\psi_{0}(x)$ satisfy the hypothese of Theorem~3.1, then one can analogously prove that
\begin{equation*}
    \mathcal{H}_{1}\big((\varphi_{0}(k))^{s};\mathcal{X}\big)
    =\mathcal{H}_{1}\big((\psi_{0}(k))^{s};\mathcal{X}\big)\quad(-\infty<s<\infty),
\end{equation*}
by the fact $(\mathcal{M},\mathcal{M})=(\mathcal{C},\mathcal{C})\subset(\mathcal{L}^{p},\mathcal{L}^{p})\;(1\leq p\leq\infty)$ (see Remark~2.4).\\
}
\textbf{Theorem~3.3}~~\emph{Let $a(x)$ and $b(x)$ be polynomials with the same degree $d$. Suppose $\mathcal{A}$ and $\mathcal{B}$ are operators in $\mathcal{X}$ with multiplier sequences $\{a(k)\}_{k=0}^{\infty}$ and $\{b(k)\}_{k=0}^{\infty}$ respectively and both possess $\alpha$-th power {\rm($\alpha>0$)}. If
$a(x)$ and $b(x)$ satisfy the hypotheses of Theorem~3.1 and neither of $a(x)$ and $b(x)$ have any zero points on $(0,+\infty)$, then there hold
\begin{equation*}
\mathcal{D}_{1}(\mathcal{A}^{\alpha})=\mathcal{D}_{1}(\mathcal{B}^{\alpha}),\quad K_{\mathcal{A}^{\alpha}}(f,\delta)_{\mathcal{X}}\approx K_{\mathcal{B}^{\alpha}}(f,\delta)_{\mathcal{X}}
\end{equation*}
for all $\alpha>0$ and $\delta>0$.}

\textbf{Proof.}~~The idea comes from \cite{Dai2003}. By Theorem~3.1, there holds
\begin{equation}\label{eq32}
\mathcal{D}_{1}(\mathcal{A}^{\alpha})=\mathcal{H}((a(k))^{\alpha};\mathcal{X})=\mathcal{H}((b(k))^{\alpha};\mathcal{X})=\mathcal{D}_{1}(\mathcal{B}^{\alpha}).
\end{equation}
For $g\in \mathcal{D}_{1}(\mathcal{A}^{\alpha})$, set
     $h := \sum\limits_{k=0}^{\infty}\Big(\frac{b(k)}{a(k)}\Big)^{\alpha}Y_k(\mathcal {A}^{r}g)=\sum_{k=0}^{\infty}\big(b(k)\big)^{\alpha}Y_k(g) \sim \mathcal {B}^{r}g$.
We show that $\|h\|_{\mathcal{X}}\leq C(a,b,\alpha,n_0)\|\mathcal {A}^{r}g\|_{\mathcal{X}}$. Setting $\psi(x)=\big(b(x)/a(x)\big)^{\alpha}$\;$(x\in [0,+\infty))$, it can be verified that
    $\left|\big(\psi(x)\big)^{(l+1)}\right|\leq C(a,b,\alpha,l)(1+x)^{-(l+2)}\;(x\geq1)$,
from which it follows that
\begin{eqnarray}\label{eq53}
  \left|\delta^{l+1}\mu_{k}\right| &\leq& \left|\int_0^{1}\cdots\int_0^{1}\psi^{(l+1)}(x)\bigg|_{x=k+u_1+u_2+\cdots+u_{l+1}}du_{1} \cdots du_{2}\right|\nonumber\\
   &\leq& C(a,b,\alpha,l) \frac{1}{(1+k)^{l+2}}.
\end{eqnarray}
Thus, for $l>\lambda$, one has
\begin{eqnarray*}
  \| h\|_{\mathcal{X}} \leq \sum_{k=0}^{\infty}\big|\delta^{l+1}\mu_{k}\big|{k+l \choose l} \left\|\sigma_k^{l}(\mathcal {A}^{\alpha} g)\right\|_{\mathcal{X}}
   \leq C(a,b,\alpha,l)\left\|\mathcal {A}^{\alpha} g\right\|_{\mathcal{X}},
\end{eqnarray*}
here the first inequality uses Abel transformations $(l+1)$ times, and the second one is by (\ref{eq46}) and (\ref{eq53}).
That is,
    $\left\|\mathcal {B}^{\alpha} g\right\|_{\mathcal{X}} = \left\|\sum\limits_{k=0}^{\infty}(b(k))^{\alpha}Y_{k}g\right\|_{\mathcal{X}}\leq C(a,b,\alpha,l)\left\|\mathcal {A}^{\alpha} g\right\|_{\mathcal{X}}$.
In the same way, one has
    $\left\|\mathcal {A}^{\alpha} g\right\|_{\mathcal{X}} \leq C(a,b,\alpha,l)\left\|\mathcal {B}^{\alpha} g\right\|_{\mathcal{X}}$.
So,
    $\left\|\mathcal {A}^{\alpha} g\right\|_{\mathcal{X}} \approx \left\|\mathcal {B}^{\alpha} g\right\|_{\mathcal{X}}$
for all $g\in \mathcal {X}$.
In addition, taking into account (\ref{eq32}), one obtains that for $f\in \mathcal {X}$,
\begin{eqnarray*}
  K_{\mathcal{A}^{\alpha}}(f,\delta) &=& \inf_{g_{1}\in \mathcal{D}_{1}(\mathcal{A}^{\alpha})}\Big\{\|f-g_{1}\|_{\mathcal{X}}+{\delta}^{\alpha}\|\mathcal{A}^{\alpha}g_{1}\|_{\mathcal{X}}\Big\} \\
   &\approx& \inf_{g_{2}\in \mathcal{D}_{1}(\mathcal{B}^{\alpha})}\Big\{\|f-g_{2}\|_{\mathcal{X}}+{\delta}^{\alpha}\|\mathcal{B}^{\alpha}g_{2}\|_{\mathcal{X}}\Big\}
   =K_{\mathcal{B}^{\alpha}}(f,\delta)\quad(\delta>0).
\end{eqnarray*}
This proves Theorem~3.3.\quad$\Box$
\section{Approximation for Semigroups of Contraction Operators of Class ($\mathscr{C}_0$) on the Sphere}
\subsection{Equivalence between approximation by semigroups of exponential-type multiplier operators and $K$-functional}
First, we prove the following lemma.\\
\textbf{Lemma~4.1.1}~~\emph{Let $\{T(t)|0\leq t<\infty\}$ be a semigroup of class $(\mathscr{C}_{0})$ in $\mathscr{E}(\mathcal {X})$ and also be multiplier operators with an exponential-type sequence $\{a_{t}(k)\}_{k=0}^{\infty}$ on $\mathcal {X}$, that is, there exits $\{a(k)\}_{k=0}^{\infty}$ such that
    $a_{t}(k)=e^{a(k)t}\;(k=0,1,\dots)$.
Then for $r\in \mathbb{Z}_{+}$,
\begin{equation}\label{eq92}
    \mathcal{D}(\mathcal{A}^{r})\subset\mathcal{D}_{1}(\mathcal{A}^{r})
\end{equation}
and
    $\mathcal{A}^{r}f\sim\sum\limits_{k=0}^{\infty}(a(k))^{r}Y_{k}f\;(f\in \mathcal{D}(\mathcal{A}^{r}))$.
Particularly,
    $\mathcal{D}(\mathcal{A})=\mathcal{D}_{1}(\mathcal{A})$.
Moreover, for $f\in \mathcal{D}(\mathcal{A}^{r})$ and $g\in \mathcal {X}$ such that
    $(a(k))^{r} Y_{k}f=Y_{k}g\;(k=0,1,\dots)$,
there holds
\begin{equation}\label{eq93}
    (T(t)-I)^{r}f=\int_{0}^{t}\cdots\int_{0}^{t}T(u_1+\cdots+u_r)g \:d u_{1}\cdots d u_{r}\quad a. e.,
\end{equation}
here $\mathcal {D}(\mathcal {A})$ and $\mathcal {D}_{1}(\mathcal {A})$ is defined by {\rm(\ref{eq18})} and {\rm(\ref{eq80})} respectively.}

\textbf{Proof.}~~First we prove $\mathcal{D}_{1}(\mathcal {A})\subset\mathcal{D}(\mathcal{A})$. Set $f\in \mathcal {D}_{1}(\mathcal{A})$ and $T f\in \mathcal{X}$ such that
    $Tf=\sum\limits_{k=0}^{\infty}a(k)Y_{k}f\;(f\in \mathcal{D}_{1}(\mathcal{A}))$.
For each fixed $x\in \mathbb{S}^{n}$, $Y_{k}(f)(x)$ $(k=0,1,\dots,)$ is a bounded linear functional on $\mathcal {X}$, which can commute with Bochner integral 
Then, for $k=0,1,2,\dots$,
\begin{eqnarray*}
  && Y_{k}\left(\int_{0}^{t}T(\tau)(T f)\:d\tau\right)(x)=\int_{0}^{t}Y_{k}\left(T(\tau)(T f)\right)(x)\:d\tau
   = \int_{0}^{t}e^{a(k)\tau}Y_{k}\big(T f\big)(x)\:d\tau\\
   &=& \int_{0}^{t}e^{a(k)\tau}a(k)\left(Y_{k}f\right)(x)\:d\tau
   = (e^{a(k)t}-1)\left(Y_{k}f\right)(x)=Y_{k}\left(T(t) f-f\right)(x), 
\end{eqnarray*}
hence by uniqueness theorem, we have
\begin{equation}\label{eq67}
    \frac{T(t) -I}{t}f=\frac{1}{t}\int_{0}^{t}T(\tau)\big(T f\big)d\tau\quad a.e. .
\end{equation}
$\{T(t)|0\leq t<\infty\}$ is of class $(\mathscr{C}_{0})$ then by (\ref{eq16}),
\begin{eqnarray*}
    \left\|\frac{T(t) -I}{t}f-T f\right\|_\mathcal {X}
    &=& \left\|\frac{1}{t}\int_{0}^{t}\Big( T(\tau)\big(T f\big)-T f\Big)d\tau\right\|_\mathcal {X}\\
    &\leq& \sup_{0\leq\tau<t}\big\|T(\tau)\big(T f\big)-T f\big\|_\mathcal {X}\rightarrow0\quad(t\rightarrow0+).
\end{eqnarray*}
Therefore, $f\in \mathcal {D}(\mathcal {A})$ and
    $\mathcal {A}f=s-\!\!\lim\limits_{t\rightarrow0+}\displaystyle\frac{T(t) -I}{t}f=T f \sim\sum_{k=0}^{\infty}a(k)Y_{k}f\;(f\in \mathcal {D}(\mathcal {A}))$,
thus $\mathcal{D}_{1}(\mathcal {A})\subset\mathcal{D}(\mathcal{A})$.

Conversely, for $f\in\mathcal {D}(\mathcal {A}^{r})$ $(r\in \mathbb{Z}_{+})$,
\begin{eqnarray*}
  &&\left(\frac{e^{a(k)t}-1}{t}\right)^{r}\left(Y_{k}f\right)(x) = Y_{k}\left(\left(\frac{T(t)-I}{t}\right)^{r}f\right)(x)\\
  &=&Y_{k}\left(\int_{0}^{t}\cdots\int_{0}^{t}T(u_{1}+\cdots+u_{r})\mathcal {A}^{r}f d u_{1}\cdots d u_{r}\right)(x)\\
  &=&\int_{0}^{t}\cdots\int_{0}^{t}e^{a(k)(u_{1}+\cdots+u_{r})} d u_{1}\cdots d u_{r} Y_{k}\left(\mathcal {A}^{r}f \right)(x)\\
  &=&\left(\frac{e^{a(k)t}-1}{t}\right)^{r}(a(k))^{-r}Y_{k}\left(\mathcal {A}^{r}f \right)(x)\quad(k=0,1,2,\dots).
\end{eqnarray*}
where the second equality is by Proposition~1.1.6 in \cite[P. 11-12]{Buter_Berens1967}. Hence,
    $Y_{k}\left(\mathcal {A}^{r}f \right)(x) = (a(k))^{r}Y_{k}f\;(k=0,1,2,\dots)$.
Thus $f\in \mathcal {D}_{1}(\mathcal {A}^{r})$. So $\mathcal {D}(\mathcal {A}^{r})\subset \mathcal {D}_{1}(\mathcal {A}^{r})$.

To prove (\ref{eq93}), we notice that
\begin{eqnarray}\label{chap4eq1}
    &&Y_{k}\left(\int_{0}^{t}\cdots\int_{0}^{t}T(u_{1}+\cdots+u_{r})g d u_{1}\cdots d u_{r}\right)(x)\nonumber\\
    &=& \int_{0}^{t}\cdots\int_{0}^{t}e^{a(k)(u_{1}+\cdots+u_{r})}(a(k))^{r}\left(Y_{k}f \right)(x) d u_{1}\cdots d u_{r}\nonumber\\
    &=& \left(e^{a(k)t}-1\right)^{r}\left(Y_{k}f \right)(x)=Y_{k}\big((T(t)-I)^{r}f \big)(x)
\end{eqnarray}
with which uniqueness theorem for Laplace series yields the result.
This completes the proof of Lemma~4.1.1.\quad$\Box$\\
\textbf{Remark~4.1.2}~~\emph{Suppose the hypotheses of Lemma~4.1.1 are satisfied, by {\rm(\ref{eq92})},
\begin{equation*}
    K_{\mathcal{A}^{r}}(f,t)_{\mathcal{X}}\leq K_{\mathcal{A}^{r}}^{*}(f,t)_{\mathcal{X}}\quad(f\in \mathcal {X},\:t>0).
\end{equation*}}
\hspace{1.3 em}By Ditzian and Ivanov's method (see \cite[P. 73-76]{Ditzian_Ivanov1993}), we obtain the following theorem.\\ 
\textbf{Theorem~4.1.3}~~\emph{Let $\mathcal{X}$ and $\mathscr{E}(\mathcal{X})$ be defined in Section~2. Suppose that $\{T(t)|0\leq t<\infty\}$ is a strongly continuous semigroup of contraction operators of class {\rm($\mathscr{C}_0$)} in $\mathscr{E}(\mathcal{X})$ and $\mathcal{A}$ is its infinitesimal generator and also $T(t)$ is an exponential-type multiplier operator for each $t>0$ defined in Lemma~4.1.1. For $f\in \mathcal{X}$ and $t>0$, $T(t)f\in \mathcal{D}(\mathcal{A})=\mathcal{D}_{1}(\mathcal{A})$ and there exists some constant $N$ independent of $t$ and $f$ such that
\begin{equation}\label{eq22}
    t\|\mathcal{A}T(t)f\|_{\mathcal{X}}\leq N\|f\|_{\mathcal{X}}\quad(\mbox{for~all}~t>0).
\end{equation}
Then, for any $r\in \mathbb{Z}_{+}$, there holds
\begin{equation*}\label{eq23}
    \|\oplus^r T(t)f - f\|_{\mathcal{X}}\approx K_{\mathcal{A}^r}(f,t^{r})_{\mathcal{X}}.
\end{equation*}}
\hspace{1.3 em}To prove Theorem~4.1.3, we need the following remark that are not difficult to verify.\\
\textbf{Remark~4.1.4}~~\emph{Let $\{T(t)|\:0\leq t<\infty\}$ be a semigroup of operators of class {\rm($\mathscr{C}_0$)}. Then for any $f\in \mathcal{D}(\mathcal{A}^{r})$ {\rm($r\in \mathbb{Z}_{+}$)} and $t>0$, there holds $T(t)f\in \mathcal{D}(\mathcal{A}^r)$ and
    $\mathcal{A}^{r} T(t)f=T(t)\mathcal{A}^{r}f$.
If $T(t)f\in \mathcal{D}(\mathcal{A})$ for all $f\in \mathcal{X}$ and all $t>0$, then for $r\in\mathbb{Z}_{+}$, $T(t)f\in \mathcal{D}(\mathcal{A}^{r})$ for all $f\in \mathcal{X}$.}

\textbf{Proof of Theorem~4.1.3.} The inequality $\|\oplus^rT(t)f-f\|_{\mathcal{X}}\leq C K_{\mathcal{A}^r}(f,t^{r})_{\mathcal{X}}$ is not hard to obtain. We just give the proof of the converse inequality
    $\|\oplus^rT(t)f-f\|_{\mathcal{X}} \geq C(r) K_{\mathcal{A}^r}(f,t^r)$.
For any $g\in \mathcal{D}_{1}(\mathcal{A}^r)$, there exists $h\in \mathcal {X}$ such that
    $Y_{k} h= a(k)Y_{k}g\;(k=0,1,2,\dots)$,
then by Lemma~4.1.1,
\begin{equation}\label{eq24}
    (T(t)-I)^rg=\int_{0}^{t}\cdots\int_{0}^{t}\int_{0}^{t}T(u_1+u_2+\cdots+u_r)h\:du_1 du_2\cdots du_r,
\end{equation}
which is a Bochner integral on $[0,t]^r$. 
In the rest part of proof of Theorem~4.1.3, we also view $\mathcal {A}^r$ as $r$-th power of the infinitesimal generator of $\{T(t)|\:0\leq t<\infty\}$.
By the corollary of Hahn-Banach theorem, for $f\in \mathcal{D}(\mathcal{A}^{r+1})$, there exists $\varphi\in \mathcal{X}^{*}$ such that \begin{equation}\label{eq47}
    \varphi\left(f+\sum_{k=1}^{r}(-1)^rT(k t)f-(-1)^r t^r \mathcal{A}^{r}f\right)=\left\|f+\sum_{k=1}^{r}(-1)^rT(k t)f-(-1)^r t^r \mathcal{A}^{r}f\right\|_{\mathcal{X}}
\end{equation}
and $\|\varphi\|_{\mathcal{X}^*}=1$. Define $F:[0,\infty)\rightarrow\mathbb{C}$ by
    $F(x)=\varphi(T(x)f)\;(0\leq x<\infty)$,
then $F\in \mathcal{C}^{r+1}[0,\infty)$ ($\mathcal{C}^{r+1}[0,\infty)$ is the space consisting of all $r$ times continuously-differentiable functions on $[0,\infty)$). 
By induction, we obtain that
\begin{equation}\label{eq48}
    F^{(i)}(x)=\varphi\left(T(x)\mathcal{A}^{i}f\right)\quad(i=1,2,\dots,r+1).
\end{equation}
Thus, 
    $\left\|F^{(r+1)}(x)\right\|_{\mathcal{C}[0,\infty)}\leq \left\|\mathcal{A}^{r+1}f\right\|_{\mathcal{X}}$.
Therefore,
\begin{eqnarray}\label{eq49}
    \sup_{x\geq0}\left|F(x)+\sum_{k=1}^{r}(-1)^k{r\choose k}F(x+kt)-(-1)^r t^r F^{(r)}(x)\right|
    \leq \frac{r}{2}\:t^{r+1}\left\|\mathcal{A}^{r+1}f\right\|_{\mathcal{X}},
\end{eqnarray}
here we use the relation between finite differences and derivatives and the mean value theorem.
Setting $x=0$, then for $t>0$ and $r\in \mathbb{Z}_{+}$, there holds, by (\ref{eq47}), (\ref{eq48}) and (\ref{eq49}) that
\begin{eqnarray}\label{eq27}
  &&\bigg\|f+\sum_{k=1}^{r}(-1)^kT(k t)f-(-1)^r t^r \mathcal{A}^{r}f\bigg\|_{\mathcal{X}}
   \leq\frac{rt^{r+1}}{2}\left\|\mathcal{A}^{r+1}f\right\|_{\mathcal{X}}.
\end{eqnarray}
On the other hand,
\begin{eqnarray}\label{chap4eq2}
    t\|\mathcal{A}^{2}T((N+2)t)f\|_{\mathcal{X}}
    &\leq& \frac{N}{N+1}\|\mathcal{A}T(t)f\|_{\mathcal{X}},
\end{eqnarray}
here Remark~4.1.4 and (\ref{eq22}) are used.
Then, there holds
\begin{eqnarray*}
  t\|\mathcal{A}T((N+2)t)f\|_{\mathcal{X}}
  &\leq& \left\|T((N+2)t)f-T((N+2)t+t)f+t\mathcal{A}T((N+2)t)f\right\|_{\mathcal{X}}\\ &&+\|T((N+2)t)(T(t)f-f)\|_{\mathcal{X}}\\
   &\leq&  \frac{t^2}{2}\|\mathcal{A}^{2}T((N+2)t)f\|_{\mathcal{X}}+\|T(t)f-f\|_{\mathcal{X}}\\
   &\leq& \frac{N}{2(N+1)}t\|\mathcal{A}T((N+2)t)f\|_{\mathcal{X}}+\left(\frac{N^2}{2}+1\right)\|T(t)f-f\|_{\mathcal{X}},
\end{eqnarray*}
where the second inequality is because of (\ref{eq27}) in the case $r=1$ and the third is due to (\ref{chap4eq2}).
Therefore,
\begin{equation}\label{eq50}
    t\|\mathcal{A}T((N+2)t)f\|_{\mathcal{X}}\leq C(N)\|T(t)f-f\|_{\mathcal{X}},
\end{equation}
where $C(N)=\big((N+1)(N^2+2)\big)/(N+2)$.
Set $m=r(N+2)$ and
    $g=-\sum\limits_{k=1}^{r}(-1)^kT(k m t)f$.
By $T(t)f\in \mathcal{D}(\mathcal{A})$ and Remark~4.1.4, there holds $g\in \mathcal{D}(\mathcal{A}^k)$ for any $k\in\mathbb{Z}_{+}$. Then, it is deduced from (\ref{eq50}) that
\begin{eqnarray*}
  t^r\|\mathcal{A}^{r}g\|_{\mathcal{X}}
  &\leq& 2^r t^r\big\|\mathcal{A}^{r}T(m t)f\big\|_{\mathcal{X}} \\
   &=& 2^r t^{r-1}\bigg(t\left\|\mathcal{A}T((N+2)t)\big(\mathcal{A}^{r-1}T((m-N-2)t)f\big)\right\|_{\mathcal{X}}\bigg) \\
   &\leq& 2^r C(N) t^{r-1}\big\|\mathcal{A}^{r-1}T((m-N-2)t)(T(t)-I)f\big\|_{\mathcal{X}} \\
   &&\qquad\dots\\
   &\leq& 2^r (C(N))^r \big\|(T(t)-I)^r f\big\|_{\mathcal{X}}.
\end{eqnarray*}
Thus, by Remark~4.1.2, one has
  $K_{\mathcal{A}^r}(f,t^{r})_{\mathcal{X}}\leq K_{\mathcal{A}^r}^{*}(f,t^{r})_{\mathcal{X}} \leq \|f-g\|_{\mathcal{X}}+t^r\|\mathcal{A}^{r}g\|_{\mathcal{X}}
   \leq \big(m^{r}+(2C(N))^r\big)\|(T(t)-I)^r f\|_{\mathcal{X}}$,
which completes the proof of Theorem~4.1.3.\quad$\Box$
\subsection{Approximation by operators with exponential-type multiplier sequences}
\textbf{Definition~4.2.1}~~\emph{Let $p(x)$ be a polynomial from $\mathbb{R}$ to $\mathbb{R}$ and $0<\gamma\leq1$, the exponential-type multiplier operator on $\mathcal {X}$ with $p(x)$ and $\gamma$ defined by
\begin{equation}\label{eq62}
    T_{p}^{\gamma}(t)f:=\sum_{k=0}^{\infty}e^{-(p(k))^{\gamma}t}Y_{k}f\quad(f\in \mathcal {X})
\end{equation}
and
\begin{equation}\label{eq82}
    T_{p}^{\gamma}(0)f := f,
\end{equation}
is called regular if the coefficient of first term is positive, $p(0)=0$ and the degree of $p(x)$ is larger than $0$.}\\
\textbf{Remark~4.2.2}~~\emph{For $f\in\mathcal {X}$,
 $T_{p}^{\gamma}(t)f=f*\varphi_{p,t}^{\gamma}$,
here
\begin{equation*}\label{eq73}
\varphi_{p,t}^{\gamma}(\cos\theta)=\frac{1}{|\mathbb{S}^{n}|}
\sum_{k=0}^{\infty}e^{-(p(k))^{\gamma}t}\frac{k+\lambda}{\lambda}P_{k}^{\lambda}(\cos\theta).
\end{equation*}
Then, $\varphi_{p,t}^{\gamma}(\cos\theta)\in \mathcal{L}^{1}_{\lambda}$ and $T_{p}^{\gamma}(t)\in \mathscr{E}(\mathcal {X})$. For $r\in \mathbb{Z}_{+}$,
\begin{equation}\label{eq63}
    (\mathcal {A}_{p}^{\gamma})^{r}f \sim \sum_{k=0}^{\infty}(-(p(x))^{\gamma})^r Y_{k}f\quad\left(\left(\mathcal{A}_{p}^{\gamma}\right)^{r}f\in \mathcal {X}\right).
\end{equation}}
\textbf{Theorem~4.2.3}~~\emph{Let $\{T_{p}^{\gamma}(t)|0\leq t<\infty\}$ defined by {\rm(\ref{eq62})} and {\rm(\ref{eq82})} be regular exponential-type multiplier operators on $\mathcal {X}$ with $p(x)$ and $\gamma$. For $t\geq0$, the kernel $\varphi_{p,t}^{\gamma}(\cos\theta)$ of $T_{p}^{\gamma}(t)$ is positive.
Then $\{T_{p}^{\gamma}(t)|0\leq t<\infty\}$ forms a strongly continuous semigroup of contraction operators of class $(\mathscr{C}_{0})$ and for $t>0$ and $f\in \mathcal {X}$,
    $T_{p}^{\gamma}(t)f\in \mathcal {D}(\mathcal {A})$.
Moreover,
\begin{equation}\label{eq65}
    \big\|\mathcal {A}_p^{\gamma}\:T_{p}^{\gamma}(t)f\big\|_{\mathcal {X}}\leq \frac{N}{t}\|f\|_{\mathcal {X}},
\end{equation}
here $N$ is a constant depending only upon $n$, $\gamma$, $p(x)$ and $\mathcal {X}$.}

\textbf{Proof.}~~For $t_{1}, t_{2}>0$ and $f\in \mathbb{Z}_{+}$, one has that
\begin{eqnarray}\label{eq83}
   T_{p}^{\gamma}(t_{1})\circ T_{p}^{\gamma}(t_{2})f 
   = T_{p}^{\gamma}(t_{1}+t_{2})f,
\end{eqnarray}
and
     $\big\|\varphi_{p,t}^{\gamma}(\cos(\cdot))\big\|_{\mathcal{L}^1_{\lambda}}= \left|\mathbb{S}^{n-1}\right|\displaystyle\int_{0}^{\pi}\varphi_{p,t}^{\gamma}(\cos\theta)(\sin\theta)^{2\lambda}\:d\theta
   =1$,
which is by (\ref{eq74}), the positivity of $\varphi_{p,t}^{\gamma}(\cos\theta)$ and (\ref{eq33}).
Thus, 
\begin{equation}\label{eq66}
    \big\|T_{p}^{\gamma}(t)f\big\|_{\mathcal {X}}\leq \big\|\varphi_{p,t}^{\gamma}(\cos(\cdot))\big\|_{\mathcal{L}^1_{\lambda}}\|f\|_{\mathcal {X}}=\|f\|_{\mathcal {X}}.
\end{equation}
and also,
\begin{equation}\label{eq42}
    \lim_{t\rightarrow0+}\|T_{p}^{\gamma}(t)f-f\|_{\mathcal{X}}=0\quad(\mbox{for~ all}~ f\in \mathcal {X}),
\end{equation}
which is by (\ref{eq66}), the contraction of $T_{p}^{\gamma}(t)$ and Banach-Steinhaus theorem as well as
the fact that the collection of all spherical polynomials is dense in $\mathcal{X}$. 
By Lemma~4.1.1, 
there holds for any $f\in \mathcal {X}$ and $t>0$,
    $T_{p}^{\gamma}(t)f\in \mathcal {D}_{1}(\mathcal {A}_{p}^{\gamma})=\mathcal {D}(\mathcal {A}_{p}^{\gamma})$,
and
\begin{equation}\label{eq69}
    \mathcal{A}_{p}^{\gamma}T_{p}^{\gamma}(t)f=\sum_{k=0}^{\infty}-(p(k))^{\gamma}Y_{k}\left(T_{p}^{\gamma}(t)f\right)
    =\sum_{k=0}^{\infty}-(p(k))^{\gamma}e^{-(p(k))^{\gamma}t} Y_{k}f\quad a. e. .
\end{equation}
Then we conclude by (\ref{eq82}), (\ref{eq83}), (\ref{eq66}) and (\ref{eq42}) that $T_{p}^{\gamma}(t)$ forms a strongly continuous semigroup of contraction operators of class $(\mathscr{C}_{0})$.

Now we go to prove (\ref{eq65}).
There exists constants $c$ and $c'$ such that
\begin{equation}\label{eq70}
    cx^{\beta}\leq(p(x))^{\gamma}\leq c' x^{\beta}\quad(0<x<\infty,\;\beta=d\gamma),
\end{equation}
where $d$ is the degree of $p(x)$. 
Then,
\begin{eqnarray}\label{eq72}
  \big\|\mathcal {A}_{p}^{\gamma}T_{p}^{\gamma}(t)f\big\|_{\mathcal {X}}
  &=& \bigg\|\sum_{k=1}^{\infty}\delta^{l+1} \big((p(k))^{\gamma}e^{-(p(k))^{\gamma}t}\big) A_{k}^{l}\sigma_{k}^{l}f\bigg\|_{\mathcal {X}}\nonumber\\
  &\leq& C \sum_{k=1}^{\infty}\bigg|\delta^{l+1} \big((p(k))^{\gamma}e^{-(p(k))^{\gamma}t}\big) k^{l}\bigg|\big\|f\big\|_{\mathcal {X}},
\end{eqnarray}
where the first equality uses $p(0)=0$ 
and Abel transformations $(l+1)$ times and $l$ is a positive integer larger than $\lambda=(n-2)/2$. 

It is necessary to estimate $\left|\sum_{k=1}^{\infty}\delta^{l+1} \big((p(k))^{\gamma}e^{-(p(k))^{\gamma}t}\big) k^{l}\right|$ this moment. One can verify by induction that 
\begin{eqnarray*}
    \left(\frac{d}{dx}\right)^{l}\Big (\big(p(x)\big)^{\gamma}e^{-(p(x))^{\gamma}t}\Big)
    &=& \sum_{i=0}^{l}e^{-(p(x))^{\gamma}t}\sum_{v=1}^{N'_{l-i}}\sum_{j=1}^{N_{i}} t^{s_{iv}} \big(p(x)\big)^{(s_{iv}+1)\gamma-(m_{iv}+r_{ij})}\\
    &&\times Q_{ivj}^{d(m_{iv}+r_{ij})-(n_{iv}+i)}(x),
\end{eqnarray*}
where $0\leq r_{ij}\leq i$, $0\leq s_{iv}, m_{iv}\leq l-i$, $n_{iv}\geq l-i$ and $N_{i}$,  $N'_{i}$ are all positive integers, $Q_{ivj}^{d}$ $(d=0,1,2,\dots)$ is a polynomial with degree $d$ and $d(m_{iv}+r_{ij})-(n_{iv}+i)\geq0$.

Thus, for $x\geq1$, one has
\begin{eqnarray*}
 &&\left|\left(\frac{d}{dx}\right)^{l}\Big (\big(p(x)\big)^{\gamma}e^{-(p(x))^{\gamma}t}\Big)\right|
 \leq \sum_{i=0}^{l} e^{-(p(x))^{\gamma}t} \sum_{v=1}^{N'_{l-i}}\sum_{j=1}^{N_{i}}C_{ivj}t^{s_{iv}}x^{(s_{iv}+1)\beta-(n_{iv}+i)}\\
 &\leq& \sum_{i=0}^{l} e^{-cx^{\beta}t} \sum_{v=1}^{N'_{l-i}}C_{iv}t^{s_{iv}}x^{(s_{iv}+1)\beta-l}
 =\sum_{i=0}^{l}C_{i} t^{i} x^{(i+1)\beta-l} e^{-cx^{\beta}t}.
\end{eqnarray*}
Therefore,
\begin{eqnarray*}
  &&\left|\delta^{l+1} \big((p(k))^{\gamma}e^{-(p(k))^{\gamma}t}\big) \right|\\
  &\leq& \left|\int_{0}^{1}\cdots\int_{0}^{1}\left(\frac{d}{dx}\right)^{l+1}\Big (\big(p(x)\big)^{\gamma}e^{-(p(x))^{\gamma}t}\Big)\Big|_{x=k+u_1+u_2+\cdots+u_{l+1}}du_1\cdots du_{l+1} \right|\\
  &\leq& \sum_{i=0}^{l+1}C'_{i}t^{i}k^{(i+1)\beta-(l+1)}e^{-ck^{\beta}t},
\end{eqnarray*}
where $C'_{i}$ $(i=0,1,\cdots,l+1)$ are positive constants depending only upon $p(x)$, $\gamma$, $i$ and $l$.
Hence,
\begin{equation}\label{eq71}
\left|\sum_{k=1}^{\infty}\delta^{l+1} \big((p(k))^{\gamma}e^{-(p(k))^{\gamma}t}\big) k^{l}\right|\leq \sum_{i=0}^{l+1}C'_{i}t^{i}\sum_{k=1}^{\infty} k^{(i+1)\beta-1}e^{-ck^{\beta}t}.
\end{equation}

Now consider function $a(x)=x^{(i+1)\beta-1}e^{-cx^{\beta}t}$ $(i=0,1,\dots,m)$,
    $\displaystyle\frac{d}{dx}\big(a(x)\big)=\Big(\big((i+1)\beta-1\big)+(-c\beta t)x^{\beta}\Big)x^{(i+1)\beta-2}e^{-cx^{\beta}t}$,
thus there exists integer $k_{i}\geq0$ (may depend on $t$) such that $a(k+1)\geq a(x)\geq a(k)$ ($1\leq k\leq x\leq k+1\leq k_{i}$) and $a(k+1)\leq a(x)\leq a(k)$ ($k+1\geq x\geq k> k_{i}$), here $k$ is a positive integer.
Then from (\ref{eq71}),
\begin{eqnarray*}
    &&\left|\sum_{k=1}^{\infty}\delta^{l+1} \big(k^{\beta}e^{-ck^{\beta}t}\big) k^{l}\right|
    \leq \sum_{i=0}^{l+1}C'_{i}t^{i}\sum_{k=1}^{\infty} k^{(i+1)\beta-1}e^{-ck^{\beta}t}\\
    &\leq& \sum_{i=0}^{l+1}C'_{i}t^{i}\bigg(\sum_{k=1}^{k_{i}}\int_{k}^{k+1}e^{-cx^{\beta}t}x^{(i+1)\beta-1}dx + \sum_{k=k_{i}+1}^{\infty}\int_{k-1}^{k}e^{-cx^{\beta}t}x^{(i+1)\beta-1}dx\bigg)\\
    &\leq& \sum_{i=0}^{l+1}\Big(2C'_{i}t^{i}\Big)\int_{0}^{\infty}e^{-cx^{\beta}t}x^{(i+1)\beta-1}dx
    = \frac{N_{1}}{t},
\end{eqnarray*}
where
    $N_{1}=\left(\sum\limits_{i=0}^{l+1}2C'_{i}c^{-(i+1)}i!\right)\beta^{-1}$.
Therefore, by (\ref{eq72}), we obtain that
    $\big\|\mathcal {A}_{p}^{\gamma}T_{p}^{\gamma}(t)f\big\|_{\mathcal {X}} \leq \displaystyle\frac{N}{t}\|f\|_{\mathcal {X}}$,
here $N=C N_{1}=N(p(x),\gamma,n,\mathcal{X})$. 
The proof of Theorem~4.2.3 is completed.\quad$\Box$\\
\textbf{Theorem~4.2.4}~~\emph{Let $T_{p}^{\gamma}(t)$ $(t\geq0)$ be regular exponential-type multiplier operators on $\mathcal {X}$ with $p(x)$ and $0<\gamma<\infty$ and their kernels are positive, then there holds for $r\in \mathbb{Z}_{+}$ that
\begin{equation}\label{eq84}
    \|\oplus^r T_{p}^{\gamma}(t)f - f\|_{\mathcal{X}}\approx K_{(\mathcal{A}_{p}^{\gamma})^r}(f,t^{r})_{\mathcal{X}}
\end{equation}
for all $f\in \mathcal {X}$, here $\mathcal {A}_{p}^{\gamma}$ with multiplier operators is the infinitesimal generator of $\{T_{p}^{\gamma}(t)|0\leq t<\infty\}$ (see {\rm(\ref{eq63})}). Moreover, $\{\oplus^{r}T_{p}^{\gamma}(t)|0\leq t<\infty\}$ is saturated with order $\mathcal {O}(t^{r})$ and its saturation class is $\mathcal{H}_{1}\big(-(p(k))^{r\gamma};\mathcal{X}\big)$.}

\textbf{Proof.}~~(\ref{eq84}) follows from Theorem~4.2.3 and Theorem~4.1.3.
We now discuss the saturation property of $\{\oplus^{r}T_{p}^{\gamma}(t)|0\leq t<\infty\}$. Let $\mathcal {A}_{p}^{\gamma}$ be the infinitesimal generator of $T_{p}^{\gamma}(t)$ and $\big\{\widehat{\varphi}_{r,p,t}^{\gamma}(k)\big\}_{k=0}^{\infty}$ be the Gegenbauer coefficients of the kernel $\varphi_{r,p,t}^{\gamma}(\cos\theta)$ of $\oplus^{r}T_{p}^{\gamma}(t)$.
Then,
    $\varphi_{r,p,t}^{\gamma}(\cos\theta)=\big(1/|\mathbb{S}^{n}|\big)\sum\limits_{k=0}^{\infty}\Big(1-\big(1-e^{-(p(k))^{\gamma}t}\big)^{r}\Big)
    \displaystyle\frac{k+\lambda}{\lambda}P_{k}^{\lambda}(\cos\theta)\in \mathcal {C}_{\lambda}(\mathbb{S}^{n})$,
and
    $\widehat{\varphi}_{r,p,t}^{\gamma}(k)=((k+\lambda)/\lambda)\Big(1-\big(1-e^{-(p(k))^{\gamma}t}\big)^{r}\Big)\;(k=0,1,2,\dots)$.
Thus,
\begin{equation}\label{eq85}
    \lim_{t\rightarrow0+}\frac{\displaystyle\frac{\lambda}{k+\lambda}\:\widehat{\varphi}_{r,p,t}^{\gamma}(k)-1}{t^{r}}=-(p(k))^{r\gamma}\quad(k=0,1,2,\dots).
\end{equation}

In addition,
\begin{equation}\label{eq91}
    \left|\mathbb{S}^{n-1}\right|\int_{0}^{\pi}\varphi_{r,p,t}^{\gamma}(\cos\theta)\left(\sin\theta\right)^{2\lambda} \:d\theta=1.
\end{equation}
and
\begin{equation}\label{eq87}
    \|\oplus^{r}T_{p}^{\gamma}(t)f\|_{\mathcal{X}}\leq 2^{r}\|f\|_{\mathcal{X}}.
\end{equation}

Using Theorem~3.1 of \cite[P. 220]{Berens_Butzer_Pawelke1968} (the $c(0,\lambda)$ there is actually $|\mathbb{S}^{n-1}|/|\mathbb{S}^{n}|$ in (\ref{eq91})), by (\ref{eq85}), (\ref{eq91}) and (\ref{eq87}), one obtains that
$f\in \mathcal{H}_{1}\big(-(p(k))^{r\gamma};\mathcal{X}\big)$ if $\|\oplus^{r}T_{p}^{\gamma}(t)f-f\|_{\mathcal{X}}=\mathcal {O}(t^{r})$ and $f$ is a constant if $\|\oplus^{r}T_{p}^{\gamma}(t)f-f\|_{\mathcal{X}}=o(t^{r})$.

Conversely, suppose $f\in \mathcal{H}_{1}\big(-(p(k))^{r\gamma};\mathcal{X}\big)$.
First, for the case of $\mathcal{X}=\mathcal{L}^{p}(\mathbb{S}^{n})$, there exists $g\in \mathcal{L}^{p}(\mathbb{S}^{n})$ such that
    $-(p(x))^{r\gamma}Y_{k}f=Y_{k}g\;(k=0,1,2,\dots)$.
By Lemma~4.1.1,
  $\left\|\oplus^{r}T_{p}^{\gamma}(t)f-f\right\|_{p} = \left\|\left(T_{p}^{\gamma}(t)-I\right)^{r}f\right\|_{p}
   = \left\|\displaystyle\int_{0}^{t}\cdots\int_{0}^{t}T(u_{1}+u_{2}+\cdots+u_{r})g\:du_{1}\cdots du_{r}\right\|_{p}\nonumber
   $\\$\leq \|g\|_{p}\: t^{r}=\mathcal {O}(t^{r})$.
The proofs  for the cases of $\mathcal{X}=\mathcal{L}^{1}(\mathbb{S}^{n})$ and $\mathcal{C}(\mathbb{S}^{n})$ are similar. So we just take the case of $\mathcal{X}=\mathcal{L}^{1}(\mathbb{S}^{n})$ for example (the framework of the proof below is from \cite[P. 229-231]{Berens_Butzer_Pawelke1968}). By hypothesis that
    $f\in \mathcal{H}_{1}\big(-(p(k))^{r\gamma};\mathcal{L}^{1}(\mathbb{S}^{n})\big)$,
there exists $\mu\in \mathcal{M}(\mathbb{S}^{n})$ such that
\begin{equation}\label{eq89}
    -(p(x))^{r\gamma}Y_{k}f=Y_{k}(d\mu)\quad(k=0,1,\dots).
\end{equation}

By Remark~2.2, the convolution
    $(\varphi_{p,t}^{\gamma}*d\mu)(x):=\displaystyle\int_{\mathbb{S}^{n}}\varphi_{p,t}^{\gamma}(x\cdot y)d\mu(y)\in \mathcal{L}^{1}(\mathbb{S}^{n})$
and
    $\|\varphi_{p,t}^{\gamma}*d\mu\|_{1}\leq\|\varphi_{p,t}^{\gamma}\|_{\mathcal{L}^{1}_{\lambda}}\|\mu\|_{\mathcal{M}}=\|\mu\|_{\mathcal{M}}$.
For given $d\mu$, $h(t)=\varphi_{p,t}^{\gamma}*d\mu$ defines a vector valued function from $(0,\infty)$ to $\mathcal{L}^{1}(\mathbb{S}^{n})$
and it can be verified that for any $\varepsilon>0$,
    $h(t)=\varphi_{p,t-\varepsilon}^{\gamma}*\left(\varphi_{p,\varepsilon}^{\gamma}*d\mu\right)
    =T_{p}^{\gamma}(t-\varepsilon)h(\varepsilon)$.
Then for $0<\varepsilon\leq t_{2}<t_{1}<\infty$, one has
  $\|h(t_{1})-h(t_{2})\|_{1} =  \left\|T_{p}^{\gamma}(t_{1}-\varepsilon)h(\varepsilon)-T_{p}^{\gamma}(t_{2}-\varepsilon)h(\varepsilon)\right\|_{1}
   \leq \left\|T_{p}^{\gamma}(t_{1}-t_{2})f-f\right\|_{1}\rightarrow0\;(t_{1}\rightarrow t_{2})$,
where (\ref{eq42}) is used.
Therefore, $h(t)$ is strongly continuous in $[\varepsilon,\infty)$ $(\varepsilon>0)$.

Now, by
    $\displaystyle\int_{\varepsilon}^{t}\|h(\tau)\|_{1}d\tau\leq \int_{\varepsilon}^{t}\|\mu\|_{\mathcal{M}}d\tau<\|\mu\|_{\mathcal{M}}\: t$,
it follows that $h(t)$ is Bochner integrable on $(0,t]$ for any $t>0$. Similar with the proof of (\ref{chap4eq1}), for $k=0,1,2,\dots$,
  $Y_{k}\left(\displaystyle\int_{0}^{t}\cdots\int_{0}^{t}\left(\varphi_{p,\:(u_{1}+u_{2}+\cdots+u_{r})}^{\gamma}*d\mu\right)\:du_{1}\cdots du_{r}\right)
   = Y_{k}\left(\oplus^{r}T_{p}^{\gamma}(t)f - f\right)$, 
hence,\\ $\oplus^{r}T_{p}^{\gamma} (t)f - f = \displaystyle\int_{0}^{t} \cdots \int_{0}^{t}\left(\varphi_{p,\:(u_{1}+u_{2}+\cdots+u_{r})}^{\gamma}*d\mu\right)\:du_{1}\cdots du_{r}$,
from which it follows that
    $\|\oplus^{r}T_{p}^{\gamma}(t)f-f\|_{1}\leq \|\mu\|_{\mathcal{M}}\:t^{r}=\mathcal {O}(t^{r})$.
This completes the proof of Theorem~4.2.4.\quad$\Box$
\section{Approximation for Generalized Spherical  Abel-Poisson and Weierstrass Operators and Their Booleans}
We now apply the results of Section~4 to two special operators, the generalized spherical Abel-Poisson operators and the generalized spherical Weierstrass operators.

The generalized spherical Abel-Poisson operators (also called generalized Abel-Poisson singular integrals) in $\mathscr{E}(\mathcal{X})$ are defined as (see \cite[P. 43-47]{Bochner1955})
\begin{equation*}
    V_{t}^{\gamma}f:=\sum_{k=0}^{\infty}\exp(-k^{\gamma}t)Y_{k}f=f*v_{t}^{\gamma}\quad(0<\gamma\leq 1,\;f\in \mathcal{X}),
\end{equation*}
where $\exp(\cdot)$ is the exponential function and $v_{t}^{\gamma}(\cos\theta)$ is the kernel that
    $v_{t}^{\gamma}(\cos\theta)=
    \big(1/|\mathbb{S}^{n}|\big)\sum\limits_{k=0}^{\infty}\exp(-k^{\gamma}t)\displaystyle\frac{k+\lambda}{\lambda}P_{k}^{\lambda}(\cos\theta)\;(0\leq\theta\leq\pi)$.
For $\gamma=1$,
set $u=e^{-t}$, one has (see \cite[P. 212-213]{Berens_Butzer_Pawelke1968})
    $V_{t}^{1}f=\sum\limits_{k=0}^{\infty}u^{k}Y_{k}f=\displaystyle\int_{\mathbb{S}^n}\frac{1}{\big|\mathbb{S}^n\big|}\frac{1-u^{2}}{(u^{2}-2u(x\cdot y)+1)^{\lambda+1}}f(y)\:d\omega_{n}(y)\;(0\leq u<1)$,
which is the classical Abel-Poisson summation on the sphere. 
 For $r\in \mathbb{Z}_{+}$, the $r$-th Boolean of $V_{t}^{\gamma}$ is
\begin{equation*}
    \oplus^{r}V_{t}^{\gamma}f=f-(I-V_{t}^{\gamma})^{r}f=\sum_{k=0}^{\infty}\left(1-\left(1-e^{-k^{\gamma}t}\right)^r\right) Y_{k}f\quad(f\in \mathcal {X}).
\end{equation*}

The generalized spherical Weierstrass operators (also called generalized spherical Weierstrass singular integrals) are given by (see \cite[P. 84]{Bochner1950})
\begin{equation*}
    W_{t}^{\kappa}f=\sum_{k=0}^{\infty}\exp\big(-(k(k+2\lambda))^{\kappa}t\big)Y_{k}f=f*w_{k}^{\gamma}\quad(0<\kappa\leq1,\;f\in \mathcal{X}),
\end{equation*}
here $w_{k}^{\gamma}$ is the kernel that
    $w_{k}^{\gamma}(\cos\theta)=\big(1/|\mathbb{S}^{n}|\big)\sum\limits_{k=0}^{\infty}\exp\big(-(k(k+2\lambda))^{\kappa}t\big)\displaystyle\frac{k+\lambda}{\lambda}P_{k}^{\lambda}(\cos\theta)\;(0\leq\theta\leq\pi)$.
For $r\in \mathbb{Z}_{+}$, the $r$-th Boolean of $W_{t}^{\kappa}$ is
\begin{equation*}
    \oplus^{r}W_{t}^{\kappa}f=f-(I-W_{t}^{\kappa})^{r}f=\sum_{k=0}^{\infty}\left(1-\left(1-e^{-(k(k+2\lambda))^{\kappa}t}\right)^r\right) Y_{k}f\quad(f\in \mathcal {X}).
\end{equation*}

The kernel of $V_{t}^{\gamma}$ and $W_{t}^{\kappa}$ are both positive, that is, for $0\leq\theta\leq\pi,\;t>0,\;0<\gamma\leq1\;0<\kappa\leq1$,
  $v_{t}^{\gamma}(\cos\theta) \geq 0$ and 
  $w_{t}^{\kappa}(\cos\theta) \geq 0$ 
which were proved in \cite{Bochner1950}, \cite[P. 43-47]{Bochner1955} and \cite{Kutter1947}.
Therefore, by Theorem~4.2.3, one has the follow lemma.\\
\textbf{Lemma~5.1}~~\emph{$V_{t}^{\gamma}$ and $W_{t}^{\kappa}$ are both regular exponential-type multiplier operators with positive kernels and with $p(x)=x$ and $p(x)=x(x+2\lambda)$ respectively and are both strongly continuous semigroups of contraction operators of class $(\mathscr{C}_{0})$ and satisfy the Bernstein-type inequality {\rm(\ref{eq65})}.}

Now we prove the equivalence between the approximation of the two operators and moduli of smoothness on the sphere.\\
\textbf{Theorem~5.3}~~\emph{Let $\{V_{t}^{\gamma}|0\leq t<\infty\}$\; \mbox{\rm($0<\gamma\leq1$)} be generalized spherical Abel-Poisson operators in $\mathscr{E}(\mathcal{X})$. Then for any $0<\gamma\leq1$ and $r\in \mathbb{Z}_{+}$, there holds for all $f\in \mathcal{X}$ that
\begin{equation*}\label{eq44}
    \left\|\oplus^{r}V_{t}^{\gamma}f-f\right\|_\mathcal{X}\approx\omega^{r\gamma}(f,t^{1/\gamma})_\mathcal{X}.
\end{equation*}}
\hspace{1.3 em}\textbf{Proof.}~~Set $\mathcal{V}^{\gamma}$ the infinitesimal generator of $\{V_{t}^{\gamma}|0\leq t<\infty\}$. By (\ref{eq63}), for $(\mathcal{V}^{\gamma})^{r}f\in\mathcal{X}$ and $r\in \mathbb{Z}_{+}$, there holds
    $(\mathcal{V}^{\gamma})^{r}f \sim \sum\limits_{k=0}^{\infty}\left(-k^{\gamma}\right)^{r}Y_{k}f$.
By Lemma~5.1 and Theorem~4.2.4,
    $\left\|\oplus^{r}V_{t}^{\gamma}f-f\right\|_{\mathcal{X}} \approx K_{(\mathcal{V}^{\gamma})^{r}}(f,t^{r})_{\mathcal{X}}$.

In Theorem~3.3, set $a(x)=(-1)^{2/\gamma}x^{2}$, $b(x)=-x(x+2\lambda)$ and $\alpha=(r\gamma)/2$, as $\lim\limits_{x\rightarrow\infty}\big((-1)^{2/\gamma}a(x)\big)/\big((-1)b(x)\big)=1$, $a(0)=b(0)=0$ and $g(t)=(1+2\lambda t)$, $(g(t))^{-1}=(1+2\lambda t)^{-1}\in C^{(2\lambda+2)}[0,+\infty)$, we obtain that
    $K_{(\mathcal{V}^{\gamma})^{r}}(f,t)_{\mathcal{X}}\approx K_{D^{\frac{r\gamma}{2}}}(f,t)_{\mathcal{X}}\;(t>0)$,
then with (\ref{eq45}) one obtains
  $\left\|\oplus^{r}V_{t}^{\gamma}f-f\right\|_{\mathcal{X}} \approx K_{(\mathcal{V}^{\gamma})^{r}}(f,t^{r})_{\mathcal{X}}
  \approx K_{D^{\frac{r\gamma}{2}}}(f,t^{r})_{\mathcal{X}} \approx {\omega}^{r\gamma}(f,t^{1/\gamma})_{\mathcal{X}}$.
The proof of Theorem~5.3 is completed.\quad$\Box$

We have the following similar result for generalized spherical Weierstrass operators.\\
\textbf{Theorem~5.4}~~\emph{Let $\{W_{t}^{\kappa}|0\leq t<\infty\}$\; \mbox{\rm($0<\kappa\leq1$)} be generalized spherical Weierstrass operators in $\mathscr{E}(\mathcal{X})$. Then for any $0<\kappa\leq1$ and $r\in \mathbb{Z}_{+}$, there holds for all $f\in \mathcal{X}$ that
\begin{equation*}\label{eq44}
    \left\|\oplus^{r}W_{t}^{\kappa}f-f\right\|_\mathcal{X}\approx\omega^{2r\kappa}(f,t^{1/{(2\kappa)}})_\mathcal{X}.
\end{equation*}}
Finally, we discuss the saturation properties for the Booleans of $V_{t}^{\gamma}$ and $W_{t}^{\kappa}$.\\
\textbf{Theorem~5.5}~~\emph{For $r\in \mathbb{Z}_{+},\;t>0$, $\oplus^{r}V_{t}^{\gamma}$ {\rm($0<\gamma\leq1$)} and $\oplus^{r}W_{t}^{\kappa}$ {\rm($0<\kappa\leq1$)}, the following statements are true.\\
(i)~~$\oplus^{r}V_{t}^{\gamma}$ is saturated with $\mathcal{O}(t^{r})$ and its saturation class is $\mathcal{H}_{1}(-k^{r\gamma};\mathcal{X})$;\\
(ii)~~$\oplus^{r}W_{t}^{\kappa}$ is saturated with $\mathcal{O}(t^{r})$ and its saturation class is $\mathcal{H}_{1}(-(k(k+2\lambda))^{r\kappa};\mathcal{X})$;\\
(iii)~~$\mathcal{H}_{1}(-k^{r\gamma};\mathcal{X})=\mathcal{H}_{1}(-(k(k+2\lambda))^{r\kappa};\mathcal{X})$ if ~$0<\gamma=2\kappa\leq1$.}

\textbf{Proof.}~~(i) and (ii) are deduced from Remark 5.1 and Theorem~4.2.4.
(iii) follows from Remark~3.2 by setting $\psi_{0}(x)=(-1)^{1/(r\kappa)}x^2$, $\varphi_{0}(x)=(-1)^{1/(r\kappa)}x(x+2\lambda)$ and $s=r\kappa$.\quad$\Box$

\end{document}